\documentclass[english]{article}
\usepackage[LGR,T1]{fontenc}
\usepackage{mathtools}
\usepackage{amsmath}
\usepackage{amsthm}
\usepackage{amssymb}
\usepackage{stackrel}
\usepackage{graphicx}

\makeatletter

\newcommand{\noun}[1]{\textsc{#1}}

\ProvideTextCommand{\~}{LGR}[1]{\char126#1}

\providecommand{\tabularnewline}{\\}

\theoremstyle{plain}
\newtheorem{thm}{\protect\theoremname}[section]
\theoremstyle{remark}
\newtheorem{rem}[thm]{\protect\remarkname}
\newtheorem{corollary}{Corollary}
\newtheorem{proposition}{Proposition}

\usepackage{bbold}
\usepackage[all]{xy}

\makeatother

\usepackage{babel}
\providecommand{\remarkname}{Remark}
\providecommand{\theoremname}{Theorem}

\begin{document}
\title{\huge{\noun{octonionic planes and real forms of $\text{G}_{2}$, $\text{F}_{4}$ and $\text{E}_{6}$}}}
\vskip 2.0 cm

\author{ D. Corradetti$^1$,  A. Marrani$^2$, D. Chester$^3$  and  R. Aschheim$^{3}$ }

\vskip 2.0 cm

\maketitle

\begin{abstract}
In this work we present a useful way to introduce the octonionic projective and hyperbolic plane $\mathbb{O}P^{2}$ through the use of Veronese
vectors. Then we focus on their relation with the exceptional Jordan algebra $\mathfrak{J}_{3}^{\mathbb{O}}$ and show that the Veronese vectors are the rank-1 elements of the algebra. We then study groups of motions over the octonionic plane recovering all real forms of
$\text{G}_{2}$, $\text{F}_{4}$ and $\text{E}_{6}$ groups and finally give a classification of all octonionic and split-octonionic planes as symmetric spaces.\\[0.2cm]
\end{abstract}

\section{Introduction}

The study of the exceptional Jordan algebra and its complexification has
been of interest in recent papers of theorethical physics. Todorov,
Dubois-Violette\cite{Todorov1} and Krasnov\cite{Krasnov} characterized
the Standard Model gauge group $\text{G}_{SM}$ as a subgroup of automorphisms
of the exceptional Jordan algebra $\mathfrak{J}_{3}\left(\mathbb{O}\right)$
while Boyle \cite{Boyle,Except} pointed to its complexification $\mathfrak{J}_{3}^{\mathbb{C}}\left(\mathbb{O}\right)$.
An equivalent well known view of the exceptional Jordan algebra is
the one of projective geometry, in which the automorphism group of $\mathfrak{J}_{3}\left(\mathbb{O}\right)$ is the group of motions of the octonionic projective plane\cite{Jordan}. Making use of
Veronese coordinates we will explore these relations and show how
all real forms of $\text{F}_{4}$ and $\text{E}_{6}$ can be recovered
as group of motions of projective or hyperbolic planes defined over
division Octonions $\mathbb{O}$ or split Octonions $\mathbb{O}_{s}$.

In sec. 2 we introduce the octonionic projective plane through
the use of Veronese coordinates, we define the projective lines and
relate the construction with the octonionic affine plane. In sec.
3 we show the correspondance with the exceptional Jordan algebra $\mathfrak{J}_{3}\left(\mathbb{O}\right)$
while in sec. 4 we show how real forms of the exceptional Lie Groups
$\text{E}_{6}$, $\text{F}_{4}$ and $\text{G}_{2}$, arise as specific groups of collineations of the octonionic projective plane. In the last section we proceed defining in a systematic way all possible octonionic planes as symmetric spaces.

\section{The Octonionic Projective Plane }

\begin{figure}
\centering{}\includegraphics[scale=0.08]{ 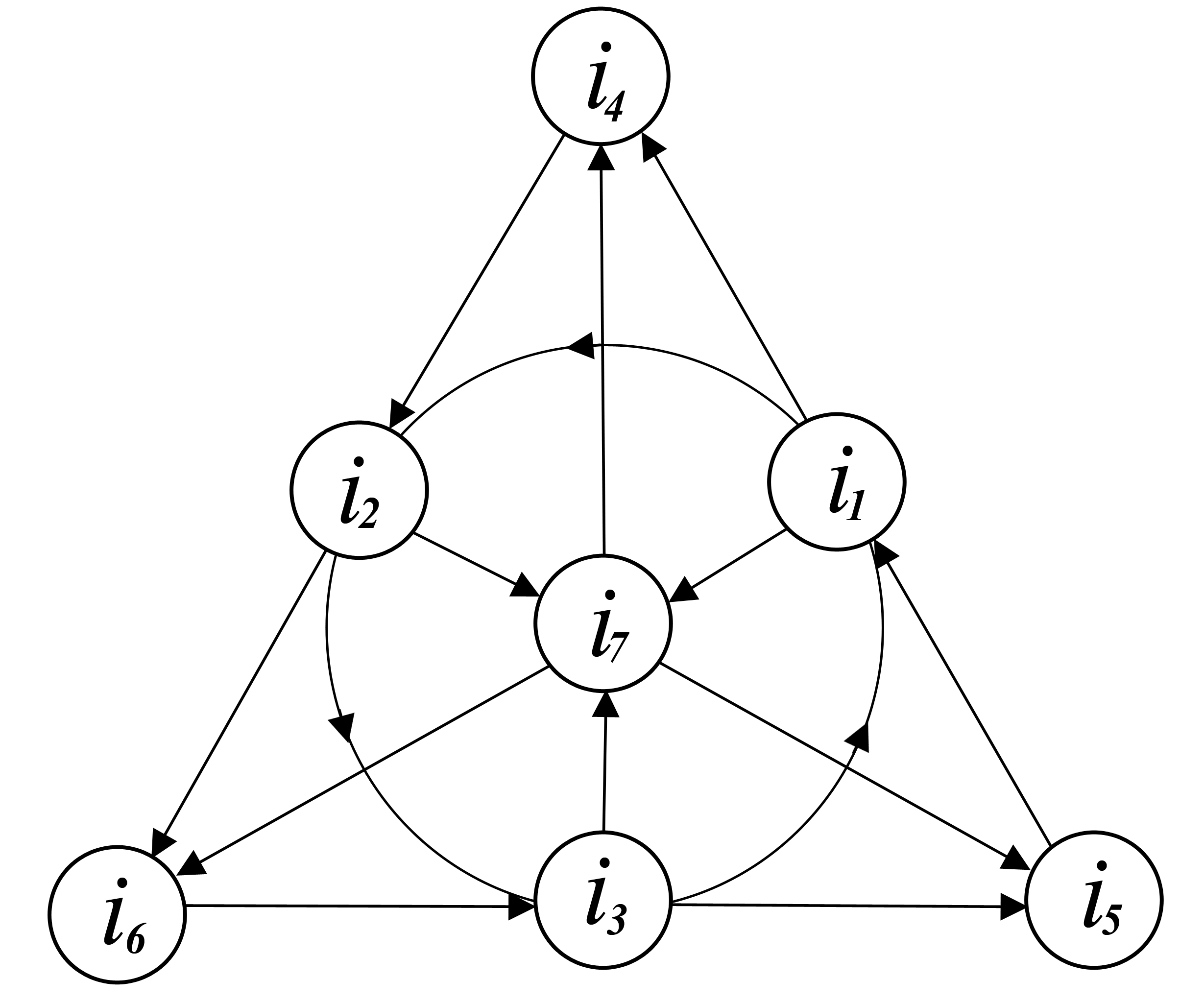}\caption{\label{fig:octonion fano plane} Multiplication rule of Octonions
$\mathbb{O}$ as real vector space $\mathbb{R}^{8}$ in the basis
$\left\{ i_{0}=1,i_{1},...,i_{7}\right\} $.}
\end{figure}
Octonions $\mathbb{O}$ are, along with Real numbers $\mathbb{R}$,
Complex numbers $\mathbb{C}$ and Quaternions $\mathbb{H}$, one of
the four \emph{Hurwitz algebras}, more specifically are the only unital
non-associative normed division algebra. A pratical way to work with
them is to consider their $\mathbb{R}^{8}$ decomposition, i.e.,
\begin{equation}
x=\sum \limits_{k=0}^7 x_{k} \textbf{i}_{k}
\end{equation}
where $\left\{ \textbf{i}_{0}=1,\textbf{i}_{1},...,\textbf{i}_{7}\right\} $ is a basis of $\mathbb{R}^{8}$
and the multiplication rules are mnemonically encoded in the \emph{Fano
plane} (Fig. \ref{fig:octonion fano plane}) along with $\textbf{i}_{k}^{2}=-1$
for $k=1,...,7$. We then define the octonionic conjugate of $x$
as 
\begin{equation}
\overline{x}\coloneqq x_{0}\textbf{i}_{0}-\sum \limits_{k=1}^7x_{k}\textbf{i}_{k}
\end{equation}
with the usual norm
\begin{equation}
\left\Vert x\right\Vert ^{2}=\overline{x}x=x_{0}^{2}+x_{1}^{2}+x_{2}^{2}+x_{3}^{2}+x_{4}^{2}+x_{5}^{2}+x_{6}^{2}+x_{7}^{2}
\end{equation}
and the inner product given by the polarisation of the norm, i.e.,
\begin{equation}
\left\langle x,y\right\rangle =\left\Vert x+y\right\Vert ^{2}-\left\Vert x\right\Vert ^{2}-\left\Vert y\right\Vert ^{2}=\overline{x}y+\overline{y}x.
\end{equation}
In respect to this norm the Octonions are a composition algebra, i.e.
$\left\Vert xy\right\Vert =\left\Vert x\right\Vert \left\Vert y\right\Vert $,
which will be of paramount importance in the following sections. Finally, we denote with $\mathbb{O}_{s}$ the split-octonionic algebra, whose definition can be found in \cite{Manogue,Rosenfeld Group2}.

\paragraph{The Projective Plane }

It is a common practice defining a projective plane over an associative
division algebra starting from a vector space over the given algebra,
e.g. $\mathbb{R}^{n+1}$, and then define the projective space as
the quotient 
\begin{equation}
\mathbb{R}P^{n}=\mathbb{R}^{n+1}/\sim
\end{equation}
where $x\sim y$ if $x$ and $y$ are multiple through the scalar
field, i.e. $\lambda x=y$, $\lambda\in\mathbb{R}$, $\lambda\neq0$,
$x,y\in\mathbb{R}^{n+1}$. But, since the algebra of Octonions is
not associative, we have that $x\left(\lambda\mu\right)\neq\left(x\lambda\right)\mu$
when $\lambda,\mu\in\mathbb{O}$. If we then try to define the equivalence
relation as above, we then might have $x\sim y=x\lambda$, and $x\sim z=x\left(\lambda\mu\right)$,
but $z$ not related to $y$ since
\begin{equation}
z=x\left(\lambda\mu\right)\neq\left(x\lambda\right)\mu=y\mu.
\end{equation}
Therefore the previous is not an equivalence relation and the quotient cannot be defined. A method for overcoming such an issue is based on determining an equivalent algebraic definition of the rank-one idempotent of the exceptional Jordan algebra in order define points in the projective plane, but here we want to use a direct and less known way to proceed making use of the Veronese vectors.

\paragraph{Veronese coordinates}

Let $V\cong\mathbb{O}^{3}\times\mathbb{R}^{3}$ be a real vector space,
with elements of the form 
\[
\left(x_{\nu};\lambda_{\nu}\right)_{\nu}=\left(x_{1},x_{2},x_{3};\lambda_{1},\lambda_{2},\lambda_{3}\right)
\]
where $x_{\nu}\in\mathbb{O}$, $\lambda_{\nu}\in\mathbb{R}$ and $\nu=1,2,3$.
A vector $w\in V$ is called \emph{Veronese} if 

\begin{align}
\lambda_{1}\overline{x}_{1} & =x_{2}x_{3},\,\,\lambda_{2}\overline{x}_{2}=x_{3}x_{1},\,\,\lambda_{3}\overline{x}_{3}=x_{1}x_{2}\label{eq:Ver-1}\\
\left\Vert x_{1}\right\Vert ^{2} & =\lambda_{2}\lambda_{3},\,\left\Vert x_{2}\right\Vert ^{2}=\lambda_{3}\lambda_{1},\left\Vert x_{3}\right\Vert ^{2}=\lambda_{1}\lambda_{2}.\label{eq:Ver-2}
\end{align}
Let $H\subset V$ be the subset of Veronese vectors. If $w=\left(x_{\nu};\lambda_{\nu}\right)_{\nu}$
is a Veronese vector then also $\mu w=\mu\left(x_{\nu};\lambda_{\nu}\right)_{\nu}$
is a Veronese vector, that is $\mathbb{R}w\subset H$. We define the
\emph{Octonionic plane} $\mathbb{O}P^{2}$ as the geometry having
this one-dimensional subspaces $\mathbb{R}w$ as points, i.e.
\begin{equation}
\mathbb{O}P^{2}=\left\{ \mathbb{R}w:w\in H\smallsetminus\left\{ 0\right\} \right\} .
\end{equation}

\begin{rem}
The point in the projective plane is defined as the equivalence class
$\mathbb{R}w$ of the Veronese vector $w$, but, in order to determine an explicit relation between points in the projective plane and rank-one idempotent
elements of the Jordan algebra $\mathfrak{J}_{3}^{\mathbb{O}}$, we
will choose when as representative of the class the vector $v=\left(y_{\nu};\xi_{\nu}\right)_{\nu}\in\mathbb{R}w$
such that $\xi_{1}+\xi_{2}+\xi_{3}=1$. Then $\left(y_{\nu};\xi_{\nu}\right)_{\nu}$ are called Veronese coordinates of the projective point.
\end{rem}

\paragraph{Projective lines}

We then define \emph{projective lines} of $\mathbb{O}P^{2}$ as the
vectors orthogonal to the points $\mathbb{R}w$. Let $\beta$ be the
bilinear form over $\mathbb{O}^{3}\times\mathbb{R}^{3}$ defined as

\begin{equation}
\beta\left(w_{1},w_{2}\right)=\sum \limits_{\nu=1}^3 \left(\left\langle x_{\nu}^{1},x_{\nu}^{2}\right\rangle +\lambda_{\nu}^{1}\lambda_{\nu}^{2}\right)
\end{equation}
where $w_{1}=\left(x_{\nu}^{1};\lambda_{\nu}^{1}\right)_{\nu}$,$w_{2}=\left(x_{\nu}^{2};\lambda_{\nu}^{2}\right)_{\nu}\in\mathbb{O}^{3}\times\mathbb{R}^{3}$.
Then, for every Veronese vector $w$, corresponding to the point $\mathbb{R}w$
in $\mathbb{O}P^{2}$, we define a line $\ell$ in $\mathbb{O}P^{2}$ as
the orthogonal space
\begin{equation}
\ell\coloneqq w^{\perp}=\left\{ z\in\mathbb{O}^{3}\times\mathbb{R}^{3}:\beta\left(z,w\right)=0\right\} .
\end{equation}

The bilinear form $\beta$ also defines the \emph{elliptic polarity},
i.e. the map $\pi^{+}$ that corresponds points to lines and lines
to points, i.e.
\begin{equation}
\pi^{+}\left(w\right)=w^{\perp},\pi^{+}\left(w^{\perp}\right)=w
\end{equation}
where the orthogonal space to a vector is defined by the bilinear
form $\beta\left(\cdot,\cdot\right)$, so that
\begin{align}
\pi^{+}: & w\longrightarrow\left\{ \beta\left(\cdot,w\right)=0\right\} \\
 & \ell\longrightarrow w
\end{align}
when $\ell$ is given by $\left\{ \beta\left(\cdot,w\right)=0\right\} $.
Explicitly, $\beta\left(w_{1},w_{2}\right)=0$ when
\begin{equation}
2\overline{x}_{1}^{1}x_{1}^{2}+2\overline{x}_{2}^{1}x_{2}^{2}+2\overline{x}_{3}^{1}x_{3}^{2}+\lambda_{1}^{1}\lambda_{1}^{2}+\lambda_{2}^{1}\lambda_{2}^{2}+\lambda_{3}^{1}\lambda_{3}^{2}=0.
\end{equation}
In addition to the elliptic polarity defined above, we then define
the \emph{hyperbolic polarity $\pi^{-}$}, which still has
\begin{equation}
\pi^{-}\left(w\right)=w^{\perp},\pi^{-}\left(w^{\perp}\right)=w
\end{equation}
but through the use of the bilinear form $\beta_{-}$ that has a change
of sign in the last coordinate, i.e. $\beta_{-}\left(w_{1},w_{2}\right)=0$
is given by
\begin{equation}
2\overline{x}_{1}^{1}x_{1}^{2}+2\overline{x}_{2}^{1}x_{2}^{2}-2\overline{x}_{3}^{1}x_{3}^{2}+\lambda_{1}^{1}\lambda_{1}^{2}+\lambda_{2}^{1}\lambda_{2}^{2}-\lambda_{3}^{1}\lambda_{3}^{2}=0.
\end{equation}
 A projective plane equipped with the hyperbolic polarity  will
be called hyperbolic plane  and denoted as $\mathbb{O}H^{2}$.
\begin{center}
\begin{figure}
\begin{centering}
\includegraphics[scale=0.3]{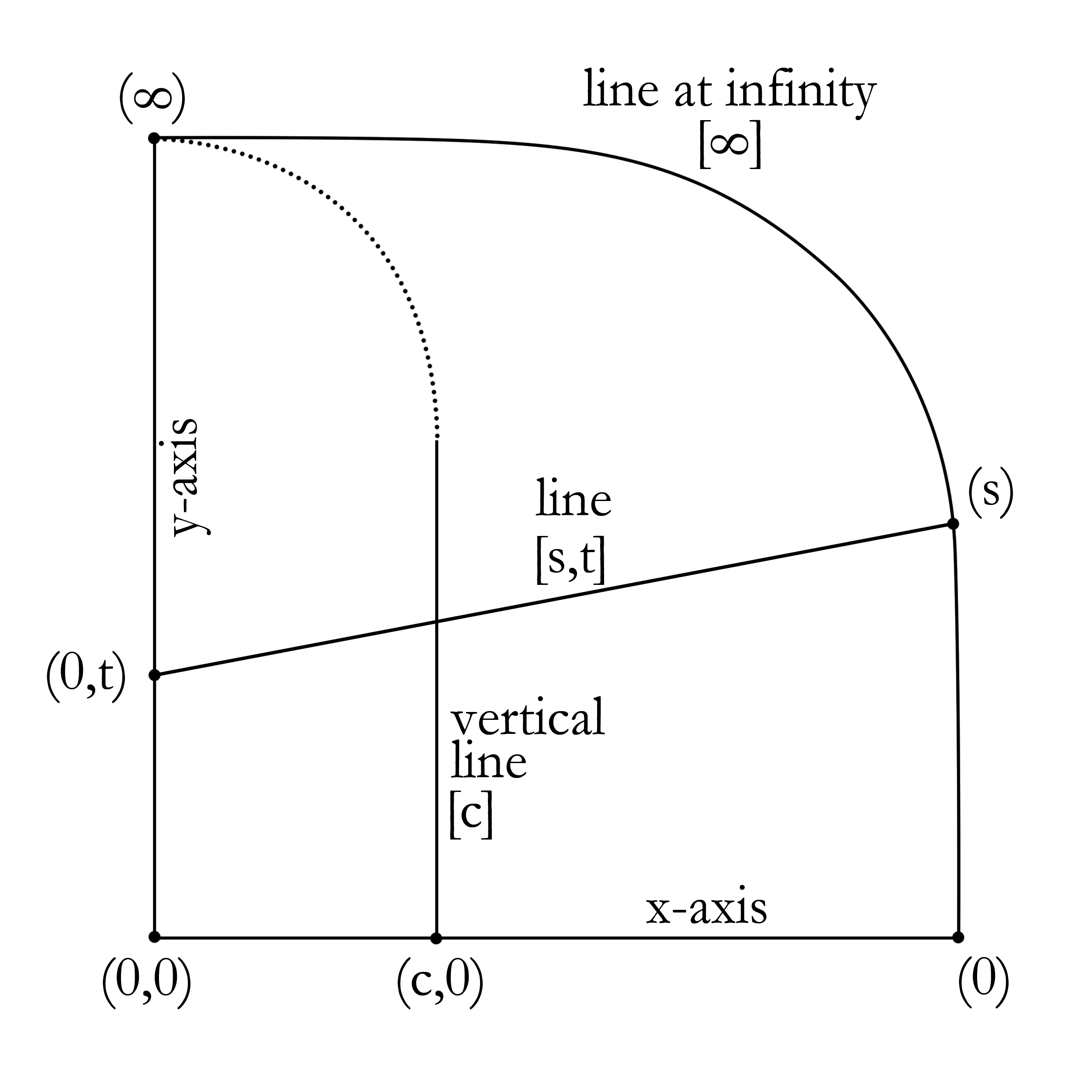}\caption{Representation of the affine plane: $\left(0,0\right)$ represents
the origin, $\left(0\right)$ the point at the infinity on the $x$-axis,
$\left(s\right)$ is the point at infinity of the line $\left[s,t\right]$
of slope $s$ while $\left(\infty\right)$ is the point at the infinity
on the $y$-axis and of vertical lines $\left[c\right]$.}
\par\end{centering}
\end{figure}
\par\end{center}

\paragraph{The Affine Plane}

The octonionic projective plane is also the completion of the octonionic
affine plane. The embedding of the affine plane can be explicited
through the use of Veronese coordinates defining the map that sends
a point $\left(x,y\right)$ of the affine plane to the projective
point $\mathbb{R}\left(x,\overline{y},y\overline{x};\left\Vert y\right\Vert ^{2},\left\Vert x\right\Vert ^{2},1\right)$,
i.e.
\begin{equation}
\left(x,y\right)\mapsto\mathbb{R}\left(x,\overline{y},y\overline{x};\left\Vert y\right\Vert ^{2},\left\Vert x\right\Vert ^{2},1\right)
\end{equation}
which is an homeomorphism. To complete the affine plane, we then have
to extend the map to another set of coordinates, i.e.
\begin{align}
\left(x\right) & \mapsto\mathbb{R}\left(0,0,x;\left\Vert x\right\Vert ^{2},1,0\right)\\
\left(\infty\right) & \mapsto\mathbb{R}\left(0,0,0;1,0,0\right).
\end{align}

\begin{rem}
To show that the above is a Veronese vector and therefore that the
map is well defined, we made essential use of alternativity of the
Octonions and fact that Octonions are a composition algebra. In case
of non-composition algebra, though the definition of the projective
and hyperbolic planes would still be valid using Veronese coordinates,
the geometry of these planes will not satisfy the basic axioms of projective
and affine geometry and therefore they would have to be considered as "generalised"
projective or hyperbolic planes.
\end{rem}

Moreover, let $\left[s,t\right]$ be a line in the affine plane $\mathbb{O}A^{2}$
of the form
\begin{equation}
\left[s,t\right]=\left\{ \left(x,sx+t\right):x\in\mathbb{O}\right\}
\end{equation}
where
$s$ is the slope of the line. Then $\left[s,t\right]$ is mapped
into the projective line orthogonal to the vector $\left(\overline{s}t,-\overline{t},-s;1,\left\Vert s\right\Vert ^{2},\left\Vert t\right\Vert ^{2}\right)$,
i.e.
\begin{equation}
\left[s,t\right]\mapsto\mathbb{R}\left(\overline{s}t,-\overline{t},-s;1,\left\Vert s\right\Vert ^{2},\left\Vert t\right\Vert ^{2}\right)^{\perp}.
\end{equation}
Vertical lines $\left[c\right]$ that are of the form $\left\{ c\right\} \times\left(\mathbb{O}\right),$
are mapped into lines of $\mathbb{O}P^{2}$ given by
\begin{equation}
\left[c\right]\mapsto\mathbb{R}\left(-c,0,0;0,1,\left\Vert c\right\Vert ^{2}\right)^{\perp}.
\end{equation}
Finally the line at infinity $\left[\infty\right]$ is mapped
to the orthogonal space of the vector 
\begin{equation}
\left[\infty\right]\mapsto\mathbb{R}\left(0,0,0;0,0,1\right)^{\perp}.
\end{equation}

\section{The Exceptional Jordan Algebra $\mathfrak{J}_{3}\left(\mathbb{O}\right)$}

The exceptional Jordan algebra $\mathfrak{J}_{3}\left(\mathbb{O}\right)$
is the algebra of Hermitian three by three octonionic matrices with
the Jordan product
\begin{equation}
X\circ Y=\frac{1}{2}\left(XY+YX\right).
\end{equation}
It is easy to see that $\mathfrak{J}_{3}\left(\mathbb{O}\right)$
is commutative, i.e. $X\circ Y=Y\circ X$ and satisfies the \emph{Jordan
identity}
\begin{equation}
\left(X^{2}\circ Y\right)\circ X=X^{2}\circ\left(Y\circ X\right).
\end{equation}
We then define the\emph{ bilinear form} 
\begin{equation}
\left(X,Y\right)=\frac{1}{2}\text{tr\ensuremath{\left(X\circ Y\right)}}
\end{equation}
the \emph{quadratic form } whose the previous bilinear form is a
polarisation
\begin{equation}
Q\left(X\right)=\frac{1}{2}\text{tr\ensuremath{\left(X^{2}\right)}}
\end{equation}
and the \emph{Freudenthal product}, i.e. 
\begin{equation}
X*Y=X\circ Y-\frac{1}{4}\left(X\text{tr\ensuremath{\left(Y\right)}}+Y\text{tr\ensuremath{\left(X\right)}}\right)+\frac{1}{4}\left(\text{tr\ensuremath{\left(X\right)}}\text{tr\ensuremath{\left(Y\right)}}-\text{tr\ensuremath{\left(X\circ Y\right)}}\right)I_{3}
\end{equation}
where $I_{3}=diag(+,+,+)$, along with the \emph{symmetric trilinear form} 
\begin{equation}
\left(X,Y,Z\right)=\frac{1}{3}\left(X,Y*Z\right)
\end{equation}
 and the\emph{ determinant}
\begin{equation}
\text{det}\left(X\right)=\frac{1}{3}\left(X,X,X\right).
\end{equation}
 Now, let $\mathbb{R}w$ be a point in the projective plane $\mathbb{O}P^{2}$,
related to a vector in $\mathbb{O}^{3}\times\mathbb{R}^{3}$ with
Veronese coordinates $w=\left(x_{\nu};\lambda_{\nu}\right)_{\nu}$
and consider the map from $V\cong\mathbb{O}^{3}\times\mathbb{R}^{3}$
into the space of three by three Hermitian matrices with octonionic
coefficients, defined as
\begin{equation}
\left(x_{\nu};\lambda_{\nu}\right)_{\nu}\mapsto\left(\begin{array}{ccc}
\lambda_{1} & x_{3} & \overline{x}_{2}\\
\overline{x}_{3} & \lambda_{2} & x_{1}\\
x_{2} & \overline{x}_{1} & \lambda_{3}
\end{array}\right).
\end{equation}
 We then have that 
\begin{equation}
\text{det}\left(X\right)=\lambda_{1}\lambda_{2}\lambda_{3}-\lambda_{1}\left\Vert x_{1}\right\Vert ^{2}-\lambda_{2}\left\Vert x_{2}\right\Vert ^{2}-\lambda_{3}\left\Vert x_{3}\right\Vert ^{2}+2\text{Re}\left(\left(x_{1}x_{2}\right)x_{3}\right)
\end{equation}
that, imposing the Veronese conditions translates to $\text{det}\left(X\right)=0$.

Moreover, let $X^{\sharp}$ be the image of a non-zero element $X$
under the adjoint ($\sharp$-)map of $\mathfrak{J}_{3}^{\mathbb{O}}$,
which is given by (cf. Example 5 of \cite{Krut})
\begin{equation}
X^{\sharp}:=\left(\begin{array}{ccc}
\lambda_{2}\lambda_{3}-\left\Vert x_{1}\right\Vert ^{2} & \overline{x_{2}}\overline{x_{1}}-\lambda_{3}x_{3} & x_{3}x_{1}-\lambda_{2}\overline{x_{2}}\\
x_{1}x_{2}-\lambda_{3}\overline{x_{3}} & \lambda_{1}\lambda_{3}-\left\Vert x_{2}\right\Vert ^{2} & \overline{x_{3}}\overline{x_{2}}-\lambda_{1}x_{1}\\
\overline{x_{1}}\overline{x_{3}}-\lambda_{2}x_{2} & x_{2}x_{3}-\lambda_{1}\overline{x_{1}} & \lambda_{1}\lambda_{2}-\left\Vert x_{3}\right\Vert ^{2}
\end{array}\right).
\end{equation}
From this explicit expression, it is immediate to realize that the
Veronese conditions are equivalent to the vanishing of $X^{\sharp}$.
Then, by the Str$\left(\mathfrak{J}_{3}\left(\mathbb{O}\right)\right)$-invariant
definition of the rank of an element of $\mathfrak{J}_{3}\left(\mathbb{O}\right)$
\cite{Jac}, one obtains that the Veronese conditions are equivalent
to the rank-1 condition for an element of $\mathfrak{J}_{3}\left(\mathbb{O}\right)$. 

Thus, from the knowledge of the orbit stratification of $\mathfrak{J}_{3}^{\mathbb{O}}$
under the non-transitive action of its reduced structure group Str$_{0}\left(\mathfrak{J}_{3}\left(\mathbb{O}\right)\right)\simeq\text{E}_{6(-26)}$,
it follows that the Veronese conditions for a non-zero element of
$\mathfrak{J}_{3}\left(\mathbb{O}\right)$ are \textit{equivalent}
to imposing that such an element belongs to the (unique) rank-1 orbit
of $\text{E}_{6(-26)}$ in $\mathfrak{J}_{3}\left(\mathbb{O}\right)$
(cf. \cite{small}, and Refs. therein).

Now we want to show that in the rank-1 (unique) orbit of $\mathfrak{J}_{3}\left(\mathbb{O}\right)$,
idempotency is equivalent to the condition of unitary trace. In order
to do that, let us consider the element $X^{2}$ and let us impose the condition
of $X$ being of rank$=1$. Since this condition is equivalent to
the Veronese conditions, one obtains
\begin{equation}
X^{2}={\scriptstyle \left(\begin{array}{ccc}
\lambda_{1}^{2}+\lambda_{1}\lambda_{2}+\lambda_{1}\lambda_{3} & \left(\lambda_{1}+\lambda_{2}+\lambda_{3}\right)x_{3} & \left(\lambda_{1}+\lambda_{2}+\lambda_{3}\right)\overline{x_{2}}\\
\left(\lambda_{1}+\lambda_{2}+\lambda_{3}\right)\overline{x_{3}} & \lambda_{2}^{2}+\lambda_{1}\lambda_{2}+\lambda_{2}\lambda_{3} & \left(\lambda_{1}+\lambda_{2}+\lambda_{3}\right)x_{1}\\
\left(\lambda_{1}+\lambda_{2}+\lambda_{3}\right)x_{2} & \left(\lambda_{1}+\lambda_{2}+\lambda_{3}\right)\overline{x_{1}} & \lambda_{3}^{2}+\lambda_{1}\lambda_{3}+\lambda_{2}\lambda_{3}
\end{array}\right)}
\end{equation}
from which it follows that $X^{2}=X$ if and only if 
\begin{equation}
\lambda_{1}+\lambda_{2}+\lambda_{3}=1
\end{equation}
i.e. if and only if $\text{tr}\left(X\right)=1$. Thus, the idempotency
condition for rank-1 elements of $\mathfrak{J}_{3}^{\mathbb{O}}$
is equivalent to the condition of unitary trace. 

\section{Lie Groups of Type $\text{G}_{2}$, $\text{F}_{4}$ and $\text{E}_{6}$ as Groups of Collineations}

We are now interested in the motions and symmetries of the octonionic projective
plane. More specifically we are interested in \emph{collineations}
that are transformations of the projective plane that send lines into
lines. If the collineation preserves the elliptic polarity or the
hyperbolic polarity is then called \emph{elliptic} or \emph{hyperbolic
motion}. Elliptic and hyperbolic motions are an equivalent characterization
of the \emph{isometries} of the projective or hyperbolic plane respectively,
thus the elliptic motion group of the projective plane  will be
indicated as $\text{Iso}\left(\mathbb{O}P^{2}\right)$ and the hyperbolic
motion group as $\text{Iso}\left(\mathbb{O}H^{2}\right)$.

\paragraph{Collineations of the Octonionic Projective Plane}

A \emph{collineation} is a bijection $\varphi$ of the set of points
of the plane onto itself, mapping lines onto lines. It is straightforward
to see that the identity map is a collineation, as the inverse $\varphi^{-1}$
and the composition $\varphi\circ\varphi'$ are if $\varphi,\varphi'$
are both collineation. Therefore the set $\text{Coll}\left(\mathbb{O}P^{2}\right)$
of collineations is a group under composition of maps. It also has a proper
subgroup of order three generated by the \emph{triality collineation}
that permutes three special points of the affine/projective octonionic
plane, i.e. the  origin of coordinate $\left(0,0\right)$,
the point at the origin of the line at infinity which has coordinate
$\left(0\right)$ and the point at infinity of the line at the infinity
which has affine coordinate $\left(\infty\right)$. In Veronese coordinates
these three points are images of the following vectors
\begin{align}
\left(0,0\right) & \longrightarrow\mathbb{R}\left(0,0,0;0,0,1\right)\\
\left(0\right) & \longrightarrow\mathbb{R}\left(0,0,0;0,1,0\right)\\
\left(\infty\right) & \longrightarrow\mathbb{R}\left(0,0,0;1,0,0\right)
\end{align}
and the triality collineation $\tau$ is given by 
\begin{equation}
\left(x_{1},x_{2},x_{3};\lambda_{1},\lambda_{2},\lambda_{3}\right)\mapsto\left(x_{2},x_{3},x_{1};\lambda_{2},\lambda_{3},\lambda_{1}\right)
\end{equation}
that is a cyclic permutation of order three that leaves invariant the Veronese vectors.
This means that it induces a bijection $\tau$ on $\mathbb{O}P^{2}$
that is unseen by the bilinear form $\beta$ and therefore maps lines
into lines, since lines are constructed as the ortogonal space of
a vector through the bilinear form $\beta$.

Let us now consider the transformations $T_{a,b}$ of $\mathbb{O}^{3}\times\mathbb{R}^{3}$
into itself defined on the Veronese coordinates as 
\begin{align}
x_{1}\longrightarrow & x_{1}+\lambda_{3}a\nonumber \\
x_{2}\longrightarrow & x_{2}+\lambda_{3}\overline{b}\nonumber \\
x_{3}\longrightarrow & x_{3}+b\overline{x_{1}}+\overline{x_{2}}\overline{a}+\lambda_{3}b\overline{a}\nonumber \\
\lambda_{1}\longrightarrow & \lambda_{1}+\left\langle \overline{x}_{2},a\right\rangle +\lambda_{3}\left\Vert b\right\Vert ^{2}\\
\lambda_{2}\longrightarrow & \lambda_{2}+\left\langle \overline{x}_{1},a\right\rangle +\lambda_{3}\left\Vert a\right\Vert ^{2}\nonumber \\
\lambda_{3}\longrightarrow & \lambda_{3}.\nonumber 
\end{align}
Those are in fact \emph{translations} on the affine plane corresponding
to the transformation $\left(x_{1},x_{2}\right)\longrightarrow\left(x_{1}+a,x_{2}+b\right)$
and they all induce collineations $T_{a,b}$ on $\mathbb{O}P^{2}$.

It can be shown that all collineations are generated by the interplay
between a translation and the conjugation of a power of the triality
collineation, i.e. are of the form 
\begin{equation}
T_{a,b},\,\,\,\tau T_{a,b}\tau^{-1},\,\,\,\tau^{2}T_{a,b}\tau^{-2}.
\end{equation}

From another perspective, collineations transform lines of $\mathbb{O}P^{2}$
in lines of $\mathbb{O}P^{2}$. This is equivalent to find all the
linear transformations $A$ of $V$ in itself
such that the image of Veronese vectors is still a Veronese vector
$A\left(H\right)\subset H$. If this condition is fulfilled, the linear
transformation $A$ in $\text{End}\left(V\right)$ will induce a
collineation $\varphi$ on $\mathbb{O}P^{2}$, i.e.
\begin{equation}
\begin{array}{ccc}
\mathbb{O}^{3}\times\mathbb{R}^{3} & \overset{A}{\longrightarrow} & \mathbb{O}^{3}\times\mathbb{R}^{3}\\
\uparrow &  & \uparrow\\
\mathbb{O}P^{2} & \overset{\varphi}{\longrightarrow} & \mathbb{O}P^{2}
\end{array}
\end{equation}
Since all linear multiple of the transformation $A$ will produce the same collineation $\varphi$, to have a bijection between linear
transformations and collineations we have to impose also $\text{det}\left(A\right)=1$.
That is that the group of collineation $\text{Coll}\left(\mathbb{O}P^{2}\right)$
is 
\begin{equation}
\text{SL}\left(V,H\right)\coloneqq\left\{ A\in\text{End}\left(V\right);A\left(H\right)\subseteq H;\text{det}\left(A\right)=1\right\} .
\end{equation}
If we also impose the preservation of the elliptic polarity, i.e.
of the bilinear form $\beta$, we will then have the group of elliptic
motion $\text{Iso}\left(\mathbb{O}P^{2}\right)$ that is
\begin{equation}
\text{SU}\left(V,H\right)=\left\{ A\in\text{End}\left(V\right);A\left(H\right)\subseteq H;\text{det}\left(A\right)=1;\text{tr}\left(A\right)=1\right\} .
\end{equation}
Those two groups are in fact two exceptional Lie Groups, i.e.
\begin{align}
\text{Coll}\left(\mathbb{O}P^{2}\right)\cong\text{SL}\left(V,H\right) & \cong E_{6}\\
\text{Iso}\left(\mathbb{O}P^{2}\right)\cong\text{SU}\left(V,H\right) & \cong F_{4}.
\end{align}
The identification of this two group is done through a direct determination
of the generators as in \cite{Manogue}; instead, we will here  follow Rosenfeld
in \cite{Boyle} focusing on the Lie algebra of the group of collineations
on $\mathbb{O}P^{2}$, i.e. $\mathfrak{coll}\left(\mathbb{O}P^{2}\right)$,
which is given by
\begin{equation}
\mathfrak{coll}\left(\mathbb{O}P^{2}\right)=\mathfrak{g}_{2}\oplus\mathfrak{a}_{3}\left(\mathbb{O}\right)
\end{equation}
where $\mathfrak{g}_{2}\cong\mathfrak{der}\left(\mathbb{O}\right)$
and $\mathfrak{a}_{3}$ are the three by three matrices on $\mathbb{O}$
with null trace, i.e. $\text{tr}\left(A\right)=0$. The dimension
count on the possible generators of this algebra, since the only condition
you have is to have null trace, i.e. $\text{tr}\left(A\right)=0$,
gives as only condition on 
\begin{equation}
A=\left(\begin{array}{ccc}
a_{1}^{1} & a_{2}^{1} & a_{3}^{1}\\
a_{1}^{2} & a_{2}^{2} & a_{3}^{2}\\
a_{1}^{3} & a_{2}^{3} & a_{3}^{3}
\end{array}\right)
\end{equation}
the condition on the trace, i.e. $a_{3}^{3}=-\left(a_{1}^{1}+a_{2}^{2}\right)$,
and therefore we have $8$ entries of dimension $8$ and $\text{dim}{}_{\mathbb{R}}\mathfrak{a}_{3}=64$.
We therefore have 
\begin{equation}
\text{dim}_{\mathbb{R}}\mathfrak{coll}\left(\mathbb{O}P^{2}\right)\cong78=64+14.
\end{equation}
Since $\mathfrak{coll}\left(\mathbb{O}P^{2}\right)$ is a Lie group,
simple and of dimension 78, then it must be of $E_{6}$ type.

\paragraph{Isometries of the Plane}
Again, following Rosenfeld \cite{Boyle} we look at the elliptic motions of $\mathbb{O}P^{2}$,
which are the collineations that preserve also the polarity $\pi^{+}$
or equivalently the form $\beta$; they are given by
\begin{equation}
\mathfrak{iso}\left(\mathbb{O}P^{2}\right)=\mathfrak{g}_{2}\oplus\mathfrak{sa}\left(3\right)
\end{equation}
where we notated $\mathfrak{sa}\left(3\right)$ the skew-Hermitian
matrices with null trace. Here the
elements of $\mathfrak{sa}\left(3\right)$ are of the form
\begin{equation}
A=\left(\begin{array}{ccc}
a_{1}^{1} & a_{2}^{1} & -\overline{a}_{3}^{1}\\
-\overline{a}_{2}^{1} & a_{2}^{2} & a_{3}^{2}\\
a_{3}^{1} & -\overline{a}_{3}^{2} & a_{3}^{3}
\end{array}\right)
\end{equation}
with $a_{ij}=\overline{a}_{ji}$, $a_{3}^{3}=-\left(a_{1}^{1}+a_{2}^{2}\right)$ and $\text{Re}\left(a_{1}^{1}\right)=\text{Re}\left(a_{2}^{2}\right)=0$.
We therefore have 3 coefficient of dimension $8$, $2$ entries of
dimension $7$ and therefore $\text{dim}_{\mathbb{R}}\mathfrak{sa}\left(3\right)=38$
so that 
\begin{equation}
\mathfrak{\text{dim}_{\mathbb{R}}\mathfrak{iso}}\left(\mathbb{O}P^{2}\right)\cong52=38+14
\end{equation}
and the group of elliptic motion $\text{Iso}\left(\mathbb{O}P^{2}\right)$, being simple and of dimension 52, is of the $F_{4}$ type.

Moreover we can proceed as in \cite{Boyle} to find the collineations
that preserve the \emph{hyperbolic polarity} $\pi^{-}$ or equivalently
the form $\beta_{-}$ we previously defined. Here the element of the
Lie algebra are of the form
\begin{equation}
A=\left(\begin{array}{ccc}
a_{1}^{1} & a_{2}^{1} & \overline{a}_{3}^{1}\\
\overline{a}_{2}^{1} & a_{2}^{2} & a_{3}^{2}\\
a_{3}^{1} & -\overline{a}_{3}^{2} & a_{3}^{3}
\end{array}\right)
\end{equation}
with $a_{1}^{1}=\left(a_{2}^{2}+a_{3}^{3}\right)$ and $\text{Re}\left(a_{1}^{1}\right)=\text{Re}\left(a_{2}^{2}\right)=0$,
therefore leading to the same count of the dimension of 
\begin{equation}
\mathfrak{\text{dim}_{\mathbb{R}}\mathfrak{iso}}\left(\mathbb{O}H^{2}\right)\cong52=38+14
\end{equation}
deducing that $\text{Iso}\left(\mathbb{O}H^{2}\right)$ is again an
$\text{F}_{4}$ type group.

\paragraph{Collineations with a Fixed Triangle or Quadrangle}

We are now interested in studying the collineations $\varphi$ on
the affine plane that fix every point of $\triangle$, i.e. $\varphi\left(\left(0,0\right)\right)=\left(0,0\right)$,
$\varphi\left(\left(0\right)\right)=\left(0\right)$ and $\varphi\left(\left(\infty\right)\right)=\left(\infty\right)$. 
\begin{proposition}
The group $\Gamma\left(\triangle,\mathbb{O}\right)$ of collineations
that fix every point of $\triangle$ are transformations of this form
\begin{align}
\left(x,y\right) & \mapsto\left(A\left(x\right),B\left(y\right)\right)\\
\left(s\right) & \mapsto\left(C\left(s\right)\right)\\
\left(\infty\right) & \mapsto\left(\infty\right)
\end{align}
where $A,B$ and $C$ are automorphisms with respect to the sum over
$\mathbb{O}$ and that satisfy
\begin{equation}
B\left(sx\right)=C\left(s\right)A\left(x\right).
\end{equation}

\begin{proof}
A collineation $\varphi$ that fixes $\left(0,0\right)$, $\left(0\right)$
and $\left(\infty\right)$, also fixes the $x$-axis and $y$-axis
and all lines that are parallel to them. This means that the first
coordinate is the image of a function that does not depend on $y$
and the second coordinate is image of a fuction that does not depend
of $x$, i.e. $\left(x,y\right)\mapsto\left(A\left(x\right),B\left(y\right)\right)$
and $\left(s\right)\mapsto\left(C\left(s\right)\right)$. Now consider
the image of a point on the line $\left[s,t\right]$. The point is
of the form $\left(x,sx+t\right)$ and its image goes to 
\begin{equation}
\left(x,sx+t\right)\mapsto\left(A\left(x\right),B\left(sx+t\right)\right).
\end{equation}
If we want this to be a collineation, the points of the line $\left[s,t\right]$
must all belong to the same line which can be easily identified setting
$x=0$, i.e. the image of $\left[s,t\right]$ is the line that joins
the points $\left(0,B\left(t\right)\right)$ and $\left(C\left(s\right)\right)$.
We now have that the condition for $\left(A\left(x\right),B\left(sx+t\right)\right)$
to be in the image of $\left[s,t\right]$ is
\begin{equation}
B\left(sx+t\right)=C\left(s\right)A\left(x\right)+B\left(t\right).\label{eq:Condition B=00003DCA}
\end{equation}
Now, if $B$ is an automorphism with respect to the sum over $\mathbb{O}$,
we then have the condition $B\left(sx\right)=C\left(s\right)A\left(x\right)$.
Conversely if $B\left(sx\right)=C\left(s\right)A\left(x\right)$ is
true that $B\left(sx+t\right)=B\left(sx\right)+B\left(t\right)$, and
$B$ is an automorphism with respect to the sum. 
\end{proof}\end{proposition}

Let us consider the quadrangle $\square$ given by the points $\left(0,0\right)$,
$\left(1,1\right)$, $\left(0\right)$ and $\left(\infty\right)$,
that is $\square=\triangle\cup\left\{ \left(1,1\right)\right\} $,
and consider the collineations that fix the $\square$. Since in
addition to the previous case we also have to impose 
\begin{equation}
\left(1,1\right)\mapsto\left(A\left(1\right),B\left(1\right)\right)=\left(1,1\right)
\end{equation}
then $C\left(1\right)=1$ and, therefore $A=B=C$ and therefore $A$
is an automorphism of $\mathbb{O}$. We then have the following 
\begin{proposition}
The collineations that fix every point of $\square$ are transformations
of the type
\begin{align}
\left(x,y\right) & \mapsto\left(A\left(x\right),A\left(y\right)\right)\\
\left(s\right) & \mapsto\left(A\left(s\right)\right)\\
\left(\infty\right) & \mapsto\left(\infty\right)
\end{align}
where $A$ is an automorphism of $\mathbb{O}$. 
\end{proposition}

Moreover, since $\text{Aut}\left(\mathbb{O}\right)=\text{G}_{2\left(-14\right)}$ and $\text{Aut}\left(\mathbb{O}_s\right)=\text{G}_{2\left(2\right)}$ \cite{Yokota}, we have the following

\begin{corollary}
The group of collineations $\Gamma\left(\square,\mathbb{O}\right)$
that fix $\left(0,0\right)$, $\left(1,1\right)$, $\left(0\right)$
and $\left(\infty\right)$ is isomorphic to $\text{Aut}\left(\mathbb{O}\right)$. Therefore  $\Gamma\left(\square,\mathbb{O}\right)$ is isomorphic to $\text{G}_{2\left(-14\right)}$,
while in the case of split octonions $\mathbb{O}_{s}$ is isomorphic
to $\text{G}_{2\left(2\right)}$.
\end{corollary}

It can be shown that the group of collineations $\Gamma\left(\triangle,\mathbb{O}\right)$
is in fact the double cover of $SO_{8}\left(\mathbb{R}\right)$, i.e.
$\text{Spin}_{8}\left(\mathbb{R}\right)$, that we define here as
\begin{equation}
\text{Spin}\left(\mathbb{O}\right)=\left\{ \left(A,B,C\right)\in O^{+}\left(\mathbb{O}\right)^{3}:A\left(xy\right)=B\left(x\right)C\left(y\right)\,\,\,\forall x,y\in\mathcal{\mathbb{O}}\right\} 
\end{equation}

where $O^{+}$ is the connected component of the orthogonal group with the identity.

\begin{proposition}
The Lie algebra $Lie\left(\Gamma\left(\triangle,\mathbb{O}\right)\right)$
of the group of collineation that fixes $\left(0,0\right),\left(0\right)$ and 
$\left(\infty\right)$ 
\begin{equation}
\mathfrak{tri}\left(\mathbb{O}\right)=\left\{ \left(T_{1},T_{2},T_{3}\right)\in\mathfrak{so}\left(\mathbb{\mathbb{O}}\right)^{3}:T_{1}\left(xy\right)=T_{2}\left(x\right)y+xT_{3}\left(y\right)\right\}
\end{equation}
 while the Lie algebra $Lie\left(\Gamma\left(\square,\mathbb{O}\right)\right)$
of the group of collineation that fixes $\left(0,0\right)$, $\left(1,1\right)$,
$\left(0\right)$ and $\left(\infty\right)$ is 

\begin{equation}
\mathfrak{der}\left(\mathbb{O}\right)=\left\{ T\in\mathfrak{so}\left(\mathbb{O}\right):T\left(xy\right)=T\left(x\right)y+xT\left(y\right)\right\} .
\end{equation}

\begin{proof}
$\Gamma\left(\triangle,\mathbb{O}\right)$ is a Lie group since it
is a closed subgroup of the Lie group of collineations. We will find
directly its Lie algebra considering the elements $A,B,C\in\Gamma\left(\triangle,\mathbb{O}\right)$
in a neighbourhood of the identity and writing them as
\[
\left(A,B,C\right)\longrightarrow\left(\text{Id}+\epsilon T_{1},\text{Id}+\epsilon T_{2},\text{Id}+\epsilon T_{3}\right)
\]
 where $T_{1},T_{2},T_{3}\in\mathfrak{so}\left(\mathbb{\mathbb{O}}\right)$.
Imposing the condition $A\left(xy\right)=B\left(x\right)C\left(y\right)$
and then we obtain 
\begin{align}
\left(\text{Id}+\epsilon T_{1}\right)\left(xy\right) & =\left(\text{Id}+\epsilon T_{2}\right)\left(x\right)\left(\text{Id}+\epsilon T_{3}\right)\left(x\right)
\end{align}
 which, considering $\epsilon^{2}=0$, yields to 
\begin{equation}
T_{1}\left(xy\right)=T_{2}\left(x\right)y+xT_{3}\left(y\right).
\end{equation}
The second part of the proposition is obtained imposing $T_{1}=T_{2}=T_{3}=T$.
\end{proof}\end{proposition}
We then have the following
\begin{align}
\Gamma\left(\triangle,\mathbb{O}\right)\cong & \text{Spin}\left(\mathbb{O}\right)\cong\text{Spin}_{8}\left(\mathbb{R}\right)\\
\Gamma\left(\square,\mathbb{O}\right)\cong & \text{Aut}\left(\mathbb{O}\right)\cong\text{G}_{2\left(-14\right)}
\end{align}
and, passing to Lie algebras, we obtain 
\begin{align}
Lie\left(\Gamma\left(\triangle,\mathbb{O}\right)\right)\cong & \mathfrak{tri}\left(\mathbb{O}\right)\cong\mathfrak{so}\left(\mathbb{\mathbb{O}}\right)\\
Lie\left(\Gamma\left(\square,\mathbb{O}\right)\right)\cong & \mathfrak{der}\left(\mathbb{O}\right)\cong\mathfrak{g}_{2(-14)}.
\end{align}

By considering the split Octonions $\mathbb{O}_s$, previous formulas yield to 
\begin{align}
\Gamma\left(\triangle,\mathbb{O}_s\right)\cong & \text{Spin}\left(\mathbb{O}_s\right)\cong\text{Spin}_{(4,4)}\left(\mathbb{R}\right)\\
\Gamma\left(\square,\mathbb{O}\right)\cong & \text{Aut}\left(\mathbb{O}\right)\cong\text{G}_{2\left(2\right)}
\end{align}
and, passing to Lie algebras, we obtain 
\begin{align}
Lie\left(\Gamma\left(\triangle,\mathbb{O}_s\right)\right)\cong & \mathfrak{tri}\left(\mathbb{O}_s\right)\cong\mathfrak{so}\left(\mathbb{\mathbb{O}_s}\right)\cong\mathfrak{so}_{4,4}\\
Lie\left(\Gamma\left(\square,\mathbb{O}\right)\right)\cong & \mathfrak{der}\left(\mathbb{O}\right)\cong\mathfrak{g}_{2(2)}.
\end{align}

Resuming all the findings, following Yokota \cite[p.105]{Yokota} in the definition of real forms of $\text{E}_{6}$, we then have the following motion groups arising from the octonionic
and split-octonionic projective and hyperbolic plane, i.e.
\begin{center}
\begin{tabular}{|c|c|c|c|}
\hline 
Proj. Space & Collineation group & Isometry group & $\Gamma\left(\square\right)$ \tabularnewline
\hline 
\hline 
$\mathbb{O}P_{\mathbb{C}}^{2}$ & $\text{E}_{6}^{\mathbb{C}}$ & $\text{F}_{4}^{\mathbb{C}}$&$G_2^{\mathbb{C}}$\tabularnewline
\hline 
$\mathbb{O}P^{2}$ & $\text{E}_{6\left(-26\right)}$ & $\text{F}_{4\left(-52\right)}$&$G_{2(-14)}$\tabularnewline
\hline 
$\mathbb{O}_{s}P^{2}$ & $\text{E}_{6\left(6\right)}$ & $\text{F}_{4\left(4\right)}$&$G_{2(2)}$\tabularnewline
\hline 
$\mathbb{O}_{s}H^{2}$ & $\text{E}_{6\left(2\right)}$ & $\text{F}_{4\left(4\right)}$&$G_{2(2)}$\tabularnewline
\hline 
$\mathbb{O}H^{2}$ & $\text{E}_{6\left(-14\right)}$ & $\text{F}_{4\left(-20\right)}$&$G_{2(-14)}$\tabularnewline
\hline 
\end{tabular}
\par\end{center}

\section{Classification of the Octonionic Projective Planes}

Thus, the space of rank-1 idempotent elements of $\mathfrak{J}_{3}\left(\mathbb{O}\right)$
enjoys the following expression as an homogeneous space 
\begin{equation}
\frac{\text{F}_{4(-52)}}{\text{Spin}_{9}},\label{OP^2}
\end{equation}
which is a compact Riemannian symmetric space, of (geodesic) rank
$=1$ and of real dimension
\begin{equation}
\text{dim}_{\mathbb{R}}\left(\frac{\text{F}_{4(-52)}}{Spin_{9}}\right)=\text{dim}_{\mathbb{R}}\left(\text{F}_{4(-52)}\right)-\text{dim}_{\mathbb{R}}\left(\text{Spin}_{9}\right)=52-36=16,
\end{equation}
as expected from the number of degrees of freedom characterizing rank-1
idempotents of $\mathfrak{J}_{3}\left(\mathbb{O}\right)$ itself.
Since the unitary trace condition is imposed on top of Veronese conditions,
the coset (\ref{OP^2}) is a (proper) submanifold of $\mathcal{O}_{\text{rank}~1}\left(\mathfrak{J}_{3}\left(\mathbb{O}\right)\right)$,
i.e.
\begin{equation}
\frac{\text{F}_{4(-52)}}{\text{Spin}_{9}}\subset\frac{\text{E}_{6(-26)}}{\text{Spin}_{9,1}\ltimes\mathbb{R}^{16}}\text{.}
\end{equation}

The space (\ref{OP^2}) of rank-1 idempotent (or, equivalently, trace-1)
elements of $\mathfrak{J}_{3}\left(\mathbb{O}\right)$ can be identified
with the (compact real form of the) octonionic projective plane $\mathbb{O}P^{2}$,
which is the largest octonionic (projective) geometry; this can also
be hinted from the fact that the tangent space to the coset $\text{F}_{4(-52)}/\text{Spin}_{9}$
transforms under the isotropy group $\text{Spin}_{9}$ as
its spinor irreducible representation $\boldsymbol{16}$, which can
indeed be realized as a pair of octonions \cite{Sudbery} :
\begin{equation}
\mathfrak{f}_{4(-52)}=\mathfrak{so}_{9}\oplus\mathbf{16}\Rightarrow T\left(\frac{\text{F}_{4(-52)}}{\text{Spin}_{9}}\right)\simeq\mathbf{16~}\text{(of~}\text{Spin}_{9}\text{)}\simeq\mathbb{O\oplus O}.
\end{equation}
Thus, one obtains that
\begin{equation}
{\textstyle \begin{array}{ccc}
\mathcal{O}_{\text{rank}~1}\left(\mathfrak{J}_{3}\left(\mathbb{O}\right)\right) & \cong & \frac{\text{E}_{6(-26)}}{\text{Spin}_{9,1}\ltimes\mathbb{R}^{16}}\\
\cup &  & \cup\\
\mathbb{O}P^{2} & \cong & \frac{\text{F}_{4(-52)}}{\text{Spin}_{9}}
\end{array}}
\end{equation}
 which gives an alternative definition of the octonionic projective plane. 
\begin{figure}
\begin{centering}
\includegraphics[scale=0.5]{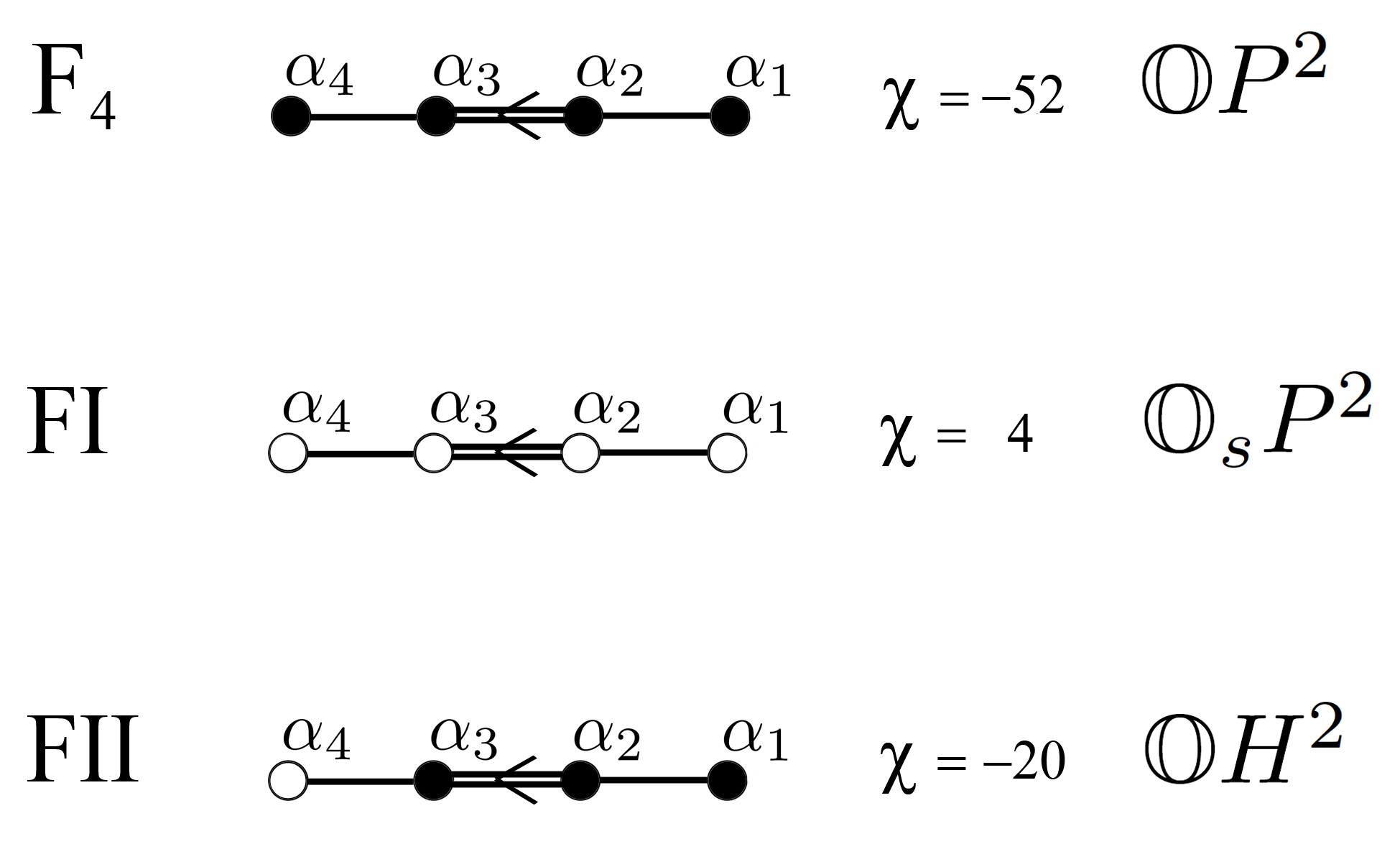}
\caption{Satake diagrams of the real forms of $\text{F}_{4}$, their character $\chi$ and the corresponding octonionic plane whose they are the isometry
group.}
\par\end{centering}
\end{figure}
From the table in the previous section, we can then classify all possible
octonionic planes. We start from the complexification of the  Cayley
plane 
\begin{eqnarray}
\mathbb{O}P^{2}\left(\mathbb{C}\right) & \simeq & \frac{\text{F}_{4}\left(\mathbb{C}\right)}{Spin_{9}\left(\mathbb{C}\right)}
\end{eqnarray}
and define three different real forms of the plane: a totally compact
real coset, that identifies with $\mathbb{O}P^{2}$; a totally non
compact which is $\mathbb{O}H^{2}$; a pseudo-Riemannian real coset
that we will define as $\mathbb{O}\widetilde{H}^{2}$. Those octonionic
planes will be defined taking as isometry group $\text{F}_{4\left(-52\right)}$
and $\text{F}_{4\left(-20\right)}$, while the last real form $\text{F}_{4\left(4\right)}$
will yield to projective planes on the split-octonionic algebra $\mathbb{O}_{s}$, i.e.
\begin{align}
\mathbb{O}P^{2} & \simeq\frac{\text{F}_{4(-52)}}{Spin_{9}}\\
\mathbb{O}H^{2} & \simeq\frac{\text{F}_{4(-20)}}{Spin_{9}}\\
\mathbb{O}\widetilde{H}^{2} & \simeq\frac{\text{F}_{4(-20)}}{Spin_{8,1}}\\
\mathbb{O}_{s}\widetilde{H}^{2} & \simeq\mathbb{O}_{s}P^{2}\simeq\mathbb{O}_{s}H^{2}\simeq\frac{\text{F}_{4(4)}}{Spin_{5,4}}.
\end{align}
Moreover, if we consider the type of the plane, i.e. the cardinality of non-compact and compact generators $\left(\#_{nc},\#_{c}\right)$, and the character $\chi$, i.e. the difference between the two,
$\chi=\#_{nc}-\#_{c}$, we then note that: the totally compact plane, i.e. the classical \emph{Cayley
plane} or the \emph{octonionic projective plane} $\mathbb{O}P^{2}$, is of type $(0,16)$ and character $\chi=16$; the totally non-compact one, i.e. \emph{hyperbolic octonionic plane} 
$\mathbb{O}H^{2}$, is of type $\left(16,0\right)$ and character
$\chi=-16$; while the other two planes named  $\mathbb{O}\widetilde{H}^{2}$ and $\mathbb{O}_{s}\widetilde{H}^{2}$ are  of type $\left(8,8\right)$ and character
$\chi=0$.

\section{Conclusions}

We have presented an explicit construction of the octonionic projective
and hyperbolic planes and showed how Lie groups of type $\text{G}_{2}$, $\text{F}_{4}$
and $\text{E}_{6}$ arise naturally as groups of motion  of such planes.
The fact that all different real forms of $\text{E}_{6}$, $\text{F}_{4}$ and $\text{G}_{2}$
can be recovered from similar constructions is of the uttermost physical
importance since different physical theories require different real
forms of Lie groups. 
Compact and non-compact forms of $\text{G}_{2}$ are notoriously known to be isomorphic to automorphisms of Octonions and split Octonions. In this paper we show how they can be thought as the subgroup of collineations that fix a quadrangle of the Projective plane over Octonions and Split-Octonions.
Different compact and non compact real forms
of $\text{E}_{6}$ and $\text{F}_{4}$ are related to different but
analogue geometric frameworks such as projective planes or hyperbolic
planes over the algebra of Octonions and split Octonions. Indeed while we recover $\text{E}_{6\left(-26\right)}$ and $\text{F}_{4\left(-52\right)}$
as collineation and isometry group of the octonionic projective plane
$\mathbb{O}P^{2}$, we have $\text{E}_{6\left(6\right)}$ and $\text{F}_{4\left(4\right)}$
for the split case $\mathbb{O}_{s}P^{2}$. The hyperbolic plane over
Octonions and split-Octonions lead to $\text{E}_{6\left(-14\right)}$
and $\text{F}_{4\left(-20\right)}$ in the octonionic case $\mathbb{O}H^{2}$
or to $\text{E}_{6\left(2\right)}$ and $\text{F}_{4\left(4\right)}$
in the split case $\mathbb{O}_{s}H^{2}$. The only real form left
out is the compact $\text{E}_{6\left(-78\right)}$ which is obtained as isometry
group of the complex Cayley plane or projective Rosenfeld plane over
Bioctonions $\left(\mathbb{C}\otimes\mathbb{O}\right)P^{2}$ \cite{CMCA}. Moreover,
we have classified all octonionic and split-octonionic projective
planes as symmetric spaces. 

\section{Acknowledgments}

The work of D.Corradetti is supported by a grant of the Quantum Gravity Research Institute. The work of AM is supported by a \textquotedblleft Maria Zambrano" distinguished researcher fellowship,
financed by the European Union within the NextGenerationEU program.

$^1$Departamento de Matemática\\ 
Universidade do Algarve\\
Campus de Gambelas\\ 
8005-139 Faro,  Portugal\\
email: a55499@ualg.pt\\[4pt]\\
$^2$Instituto de Física Teorica, Dep.to de Física,\\
Universidad de Murcia, Campus de Espinardo, E-30100, Spain\\
email:jazzphyzz@gmail.com\\[4pt]
\hspace*{2.6cm} \\
$^3$ Quantum Gravity Research,\\
Los Angeles, California, CA 90290, USA\\
davidC@quantumgravityresearch.com\\
raymond@quantumgravityresearch.com\\
\label{last}
\end{document}